\documentclass[12pt]{article}
\usepackage{amssymb,amsmath}
\usepackage[T2A]{fontenc} \usepackage[cp1251]{inputenc}
\usepackage[english]{babel}

\newtheorem{thrm}{Theorem}[section]
\newtheorem{lem}[thrm]{Lemma}
\newtheorem{prop}[thrm]{Proposition}
\newtheorem{cor}[thrm]{Corollary}

\newtheorem{remark}[thrm]{Remark}
\numberwithin{equation}{section}

\title{\textbf{The Skitovich-Darmois theorem for finite Abelian groups}}
\author{I.P. Mazur}

\begin{document}
\maketitle
\begin{abstract} Let $X$ be a finite Abelian group, $\xi_i, i=1,2,\ldots,n,n\geq 2,$ be independent random variables with values in $X$ and distributions $\mu_i$. Let $\alpha_{ij},i,j=1,2,\ldots,n,$ be automorphisms of $X$. We prove that the independence of n linear forms $L_j=\sum_{i=1}^{n}\alpha_{ij}\xi_i$ implies that all $\mu_i$ are shifts of the Haar distributions on some subgroups of the group $X$. This theorem is an analogue of the Skitovich-Darmois theorem for finite Abelian groups.
\end{abstract}

\section{Introduction}

The classical Skitovich-Darmois theorem states (\cite{Skit},\cite{Darm}): Let $\xi_i ,i=1,2,\ldots,n,n\geq2,$ be independent random variables, and $\alpha_i,\beta_i$ be nonzero numbers. Suppose that the linear forms $L_1=\alpha_1\xi_1+\cdots+\alpha_n\xi_n$ and $L_2=\beta_1\xi_1+\cdots+\beta_n\xi_n$ are independent. Then all random variables $\xi_i $ are Gaussian.

Ghurye and Olkin generalized the Skitovich-Darmois theorem to the case when $\xi_i$ are random vectors with values in $\mathbb{R}^m$, and $\alpha_i,\beta_i$ are nonsingular matrixes (\cite{Gh Ol}). They proved that the independence of the linear forms $L_1$ and $L_2$ implies that all $\xi_i$ are Gaussian vectors.

The Skitovich-Darmois theorem was generalized into various classes of locally compact Abelian groups such as finite, discrete, compact Abelian groups (see \cite{F 1992}-\cite{F 2000},\cite{FG 2000}-\cite{FG 2010}).
In the article we continue these researches  and study the Skitovich-Darmois theorem in the case when random variables take values in a finite Abelian group and the number of linear forms more than 2.

Throughout the article $X$ will denote a finite Abelian group  unless the contrary is explicitly specified. Let $Aut(X)$ be the group of automorphisms of the group $X$, $\mathbb{Z}(k)=\{0,1,2\ldots,k-1\}$ be the group of residue modulo $k$. Let $x\in X$. Denote by $E_x$ the degenerate distribution, concentrated at $x$. Let $K$ be a subgroup of $X$. Denote by $m_K$ the Haar distribution on $K$. Denote by $I(X)$ the set of shifts of such distributions, i.e. the distributions of the form $m_K\ast E_x$, where $K$ is a subgroup of $X$, $x\in X$. The distributions of the class $I(X)$ are called idempotent. Note that the idempotent distributions on a finite Abelian group can be regarded as analogues of the Gaussian distributions on real line.

Let $\xi_i, i=1,2,\ldots,n,n\geq 2,$ be independent random variables taking values in $X$ and with distributions $\mu_i$. Let $\alpha_j,\beta_j$ be automorphisms of $X$. Consider the linear forms $L_1=\alpha_1\xi_1+\cdots+\alpha_n\xi_n$ and $L_2=\beta_1\xi_1+\cdots+\beta_n\xi_n$.
The problem of the generalization of the Skitovich-Darmois theorem to the finite Abelian groups was considered first in \cite{F 1992}, where in particular it was proved that the class of groups, on which the independence of $L_1$ and $L_2$ implies that all $\mu_i$ are idempotent distributions is poor and consists of the groups of the form
\begin{equation}
\mathbb{Z}(2^{m_1})\times\cdots\times\mathbb{Z}(2^{m_l}),0\leq m_1<\cdots<m_l. \label{Lestnica}
\end{equation}

On the other hand if we consider two linear forms of two independent random variables, then the Skitovich-Darmois theorem is valid for an arbitrary finite Ableian group. Namely, the following theorem holds (\cite{F 2000}, see also \cite[p. 133]{F book}):
\begin{thrm}\label{Th n=2}
Let $\xi_1$ and $\xi_2$ be independent random variables with values in $X$ and distributions $\mu_1$ and $\mu_2$. Let $\alpha_i,\beta_i\in Aut(X),i=1,2$. If the linear forms $L_1=\alpha_1\xi_1+\alpha_2\xi_2$ and $L_2=\beta_1\xi_1+\beta_2\xi_2$ are independent, then $\mu_i \in I(X),i=1,2$.
\end{thrm}\label{Th 1}
In the paper we consider $n$ linear forms $L_j$ of $n$ independent random variables $\xi_i$ with values in a finite Abelian group. Coefficients of the forms are automorphisms of the group. We prove that the independence of $L_j$ implies that all $\xi_i$ have idempotent distributions.
This result generalizes Theorem~\ref{Th n=2} and can be considered as a natural analogue of the Skitovich-Darmois theorem for finite Abelian groups.

The main result of the article is the following theorem.

\begin{thrm}\label{Th 1}
Let $\xi_i, i=1,2,\ldots,n,n\geq 2,$ be independent random variables with values in a group $X$ and distributions $\mu_i$.
If the linear forms $L_j=\sum_{i=1}^{n}\alpha_{ij}\xi_i,$
where $\alpha_{ij}\in Aut(X),i,j=1,2,\ldots,n$, are independent, then $\mu_i\in I(X),i=1,2,\ldots,n$.
\end{thrm}
Note that the proof of Theorem~\ref{Th 1} differs from the proof of Theorem~\ref{Th n=2} for n=2 and does not use it.

Also we show that Theorem~\ref{Th 1} fails if we consider less than $n$ linear forms of $n$ random variables.

To prove the main theorem we will use some notions and results of abstract harmonic analysis (see \cite{H R}).
Let $Y=X^*$ be the character group of $X$. Since $X$ is a finite group, $Y\cong X$. The value of a character $y\in Y$ at $x\in X$ denote by $(x,y)$.
Let $\alpha:X \rightarrow X$ be a homomorphism. For each $y\in Y$ define the mapping $\tilde{\alpha}: Y \rightarrow Y$ by the equality
$(\alpha x,y)=(x,\tilde{\alpha}y)$ for all $x\in X, y \in Y$. The mapping $\tilde{\alpha}$ is a homomorphism. It is called an adjoint of $\alpha$. The identity automorphism of a group denote by $I$.
Let $B$ be a subgroup of $X$. Put
$A(Y,B)=\{y\in Y:(x,y)=1 $ for all $ x\in\ B\}.$
The set $A(Y,B)$ is called the annihilator of $B$ in $Y$ and $A(Y,B)$ is a subgroup of $Y$.

A subgroup $H$ of $X$ is called characteristic if the equality $\gamma H = H$ holds for all $\gamma \in Aut(X)$. Let $p$ be a prime number. We recall that an Abelian group is called an elementary $p$-group if every nonzero element of this group has order $p$. We note that every finite elementary $p$-group is isomorphic to a group of the form $(\mathbb{Z}(p))^m$ for some $m$. Put $X_{(p)}=\{x\in X:px=0\}$. Obviously, $X_{(p)}$ is an elementary $p$-group. Also it is obvious that $X_{(p)}$ is a characteristic subgroup of $X$.

Let $E$ be a finite-dimensional linear space and $\gamma$ be a linear operator acting on $E$. Denote by $\dim E$ the dimension of $E$ and by $Ker\gamma$ the kernel of $\gamma$. Let $\{E_i\}_{i=1}^{n}$ be a family of linear spaces. Denote by $\oplus_{i=1}^{n} E_i$ a direct sum of the linear spaces $E_i,i=1,2,\ldots,n$.

Let $\mu$ be a probability distribution on $X$. Denote by $\sigma(\mu)$ the support of $\mu$.
Put $\bar{\mu}(M)=\mu(-M)$, where $M\subset X,-M=\{-m: m\in M\}$.
The characteristic function of the distribution $\mu$ is defined by the formula: $$\hat{\mu}(y)=\sum_{x\in X}(x,y)\mu(\{x\}), \quad y\in Y.$$
If $\xi$ is a random variable with values in $X$ and distribution $\mu$, then $\hat{\mu}(y)=\mathbf{E}[(\xi,y)]$.
Put
$$F_{\mu}=\{y\in Y: \hat{\mu}(y)=1\}.$$
Then $F_{\mu}$ is a subgroup of $Y$, the inclusion $\sigma(\mu)\subset A(X,F_{\mu})$ holds, and  $\hat{\mu}(y+h)=\hat{\mu}(y)$ for all $y\in Y, h \in F_{\mu}$. If $K$ is a subgroup of $X$, then
\begin{equation}\label{HF Haara}
\hat{m}_K(y)=
\begin{cases}
1,& y\in A(Y,K), \\
0, & y\not\in A(Y,K).
\end{cases}
\end{equation}

\section{The lemmas}

To prove Theorem~\ref{Th 1} we need some lemmas. The proof of the next lemma uses standard arguments (see \cite[p. 93]{F book}).
\begin{lem}\label{Lem 1}
Let $\xi_i, i=1,2,\ldots,n,n\geq 2,$ be independent random variables with values in a group $X$ and distributions $\mu_i$.
Consider the linear forms
$L_j=\sum_{i=1}^{n}\alpha_{ij}\xi_i,j=1,2,\ldots,k,$
where $\alpha_{ij}$ are endomorphisms of $X$.
The linear forms $L_j$ are independent if and only if the following equality holds
\begin{equation}
\prod_{i=1}^{n}\hat{\mu}_i\left(\sum_{j=1}^{k}\tilde{\alpha}_{ij}u_j\right)=\prod_{i=1}^{n}\prod_{j=1}^{k}\hat{\mu}_i(\tilde{\alpha}_{ij}u_j),\quad u_j\in Y.\label{S-D obsh nk}
\end{equation}
\end{lem}

\textbf{Proof. }
We note that the linear forms $L_j,j=1,2,\ldots,k,$ are independent if and only if the equality

\begin{equation}
\textbf{E}\left[\prod_{j=1}^{k}\left(\sum_{i=1}^{n}\alpha_{ij}\xi_i,u_j\right)\right]=\prod_{j=1}^{k}\textbf{E}\left[\left(\sum_{i=1}^{n}\alpha_{ij}\xi_i,u_j\right)\right], \quad u_i \in Y  \label{lemma1}
\end{equation}
holds.
Taking in the account that the random variables $\xi_i$ are independent and that
$\hat{\mu}_i(y)=\mathbf{E}[(\xi_i,y)]$, we transform the left hand side of the equality (\ref{lemma1}) to the form
$$\textbf{E}\left[\prod_{j=1}^{k}\left(\sum_{i=1}^{n}\alpha_{ij}\xi_i,u_j\right)\right]=\textbf{E}\left[\prod_{i=1}^{n}\left(\xi_i,\sum_{j=1}^{k}\tilde{\alpha}_{ij}u_j\right)\right]=$$
$$=\prod_{i=1}^{n}\textbf{E}\left[\left(\xi_i,\sum_{j=1}^{k}\tilde{\alpha}_{ij}u_j\right)\right]=\prod_{i=1}^{n}\hat{\mu}_i\left(\sum_{j=1}^{k}\tilde{\alpha}_{ij}u_j\right).$$
Reasoning similar, we transform the right hand side of the equality  (\ref{lemma1}):
$$\prod_{i=1}^{n}\textbf{E}\left[\left(\sum_{j=1}^{k}\alpha_{ij}\xi_i,u_j\right)\right]=\prod_{i=1}^{n}\textbf{E}\left[\prod_{j=1}^{k}(\alpha_{ij}\xi_i,u_j)\right]=$$
$$=\prod_{i=1}^{n}\textbf{E}\left[\prod_{j=1}^{k}(\xi_i,\tilde{\alpha}_{ij}u_j)\right]=\prod_{i=1}^{n}\prod_{j=1}^{k}\textbf{E}\left[(\xi_i,\tilde{\alpha}_{ij}u_j)\right]=\prod_{i=1}^{n}\prod_{j=1}^{k}\hat{\mu}_i(\tilde{\alpha}_{ij}u_j).$$
$\blacksquare$

\begin{lem}\label{Lem 2}
Let $Y$ be a linear space, $\beta_{ij}$ be invertible linear operators acting on $Y$ and satisfying the conditions $\beta_{1j}=I, \beta_{i1}=I,i,j=1,2,\ldots,n$, where $I$ is the identity operator. Let $\{E_i\}_{i=1}^{n},\{F_i\}_{i=1}^{n}$ be families of finite-dimensional linear subspaces of $Y$ satisfying the conditions:

\begin{equation}
\beta_{ij}(E_j)\subset F_i,\quad i,j=1,2,\ldots,n, \label{alpha Ei sub Fj}
\end{equation}
\begin{equation}
\sum_{i=1}^{n}\dim F_i\leq \sum_{i=1}^{n} \dim E_i.\label{proobraz subset E_1 E_2 l1}
\end{equation}
Then $E_i=F_j=F,i,j=1,2,\ldots,n$, where $F$ is a linear subspace of $Y$ and $\beta_{ij}(F)=F$.
\end{lem}
\textbf{Proof. }
Put $\dim E_i=m_i,\dim F_i=k_i$. Then inequality (\ref{proobraz subset E_1 E_2 l1}) takes the form
\begin{equation}
\sum_{i=1}^{n} k_i\leq\sum_{i=1}^{n} m_i. \label{dim 0}
\end{equation}

Since $\beta_{ij}$ are invertible, we have
\begin{equation}
\dim \beta_{ij}(E_j)=m_j,\quad i,j=1,2,\ldots,n.\label{dim alpha}
\end{equation}

From (\ref{alpha Ei sub Fj}) and (\ref{dim alpha}) it follows that
\begin{equation}
m_i\leq k_j,\quad i,j=1,2,\ldots, n. \label{ner 1}
\end{equation}
From (\ref{ner 1}) we obtain that
$$\max_{1\leq i\leq n} m_i \leq \min_{1\leq j\leq n} k_j.$$
From this and (\ref{dim 0}) it follows that

\begin{equation}
\sum_{i=1}^n k_i\leq \sum_{i=1}^nm_i\leq n \min_{1\leq j\leq n} k_j. \label{ner 10}
\end{equation}
Hence, (\ref{ner 10}) implies that $k_j=k$ and (\ref{ner 10}) takes form
$$nk\leq \sum_{i=1}^nm_i\leq nk.$$
This implies that $\sum_{i=1}^n m_i=nk$. From this and $m_i\leq k,i=1,2,\ldots,n,$ it follows that $m_i=k,i=1,2,\ldots,n$. From this and from (\ref{alpha Ei sub Fj}) we obtain that
\begin{equation}
\beta_{ij}(E_j)=F_i,\quad i,j=1,2,\ldots,n. \label{alpha Ei=Fj}
\end{equation}

From (\ref{alpha Ei=Fj}) and the equalities $ \beta_{1j}=\beta_{i1}=I,i,j=1,2,\ldots,n,$ we infer
$$F_1=\beta_{1j}(E_j)=I(E_j)=E_j,$$
$$F_i=\beta_{i1}(E_1)=I(E_1)=E_1,$$
whence we have
\begin{equation}
E_i=F_j=F,i,j=1,2,\ldots,n,\label{Fi=Ei}
\end{equation}
where $F$ is a subspace of $Y$.
From (\ref{alpha Ei=Fj}) and (\ref{Fi=Ei}) it follows that $\beta_{ij}(F)=F$, $i,j=1,2,\ldots,n$.
$\blacksquare$

\begin{lem}\label{Lem 3}
Let $Y$ be a finite elementary $p$-group. Let $\hat{\mu}_i(y),i=1,2,\ldots,n$, $n\geq 2,$ be nonnegative characteristic functions on $Y$, satisfying the equation
\begin{equation}
\prod_{i=1}^{n}\hat{\mu}_i\left(\sum_{j=1}^{n}\beta_{ij}u_j\right)=\prod_{i=1}^{n}\prod_{j=1}^{n}\hat{\mu}_i(\beta_{ij}u_j),\quad u_j\in Y,\label{S-D obsh beta}
\end{equation}
where $\beta_{ij}\in Aut(Y),\beta_{1j}= \beta_{i1}=I,i,j=1,2,\ldots,n$.
Then $F_{\mu_i}=F,i=1,2,\ldots,n$, where $F$ is a subgroup of $Y$ and $\beta_{ij}(F)=F,i,j=1,2,\ldots,n$.
\end{lem}
\textbf{Proof. }
We note that $Y$ is a finite-dimensional linear space over the field $\mathbb{Z}(p)$. Then subgroups of $Y$ are subspaces of $Y$, and automorphisms acting on $Y$ are invertible linear operators.

Let $\pi$ be a map from $Y^n$ to $Y^n$ defined by the formula
\begin{equation}
\pi(u_1,u_2,\ldots,u_n)=\left(\sum_{j=1}^{n}\beta_{1j}u_j,\sum_{j=1}^{n}\beta_{2j}u_j,\ldots,\sum_{j=1}^{n}\beta_{nj}u_j\right),\label{Def pi}
\end{equation}
where $u_j \in Y$. Then $\pi$ is a linear operator. Generally, $\pi$ is not invertible.

Put $N=\pi^{-1}(\oplus_{i=1}^{n}F_{\mu_i})$.
Obviously,
\begin{equation}
\dim\oplus_{i=1}^{n} F_{\mu_i}\leq \dim N. \label{dim F leq dim N}
\end{equation}

Let $\phi_i$ be the projection on the $i$-th coordinate subspace of $Y^n$. Put $E_i=\phi_i(N)$. Then $E_i$ is a subspace of $Y$. We will show that the families of the subspaces $\{E_i\}_{i=1}^{n},\{F_{\mu_i}\}_{i=1}^{n}$ satisfy conditions (\ref{alpha Ei sub Fj}) and (\ref{proobraz subset E_1 E_2 l1}).

It is obvious that $N\subseteq (\oplus_{i=1}^{n} E_i)$. From this and (\ref{dim F leq dim N}) we obtain that
\begin{equation}
\dim \oplus_{i=1}^{n} F_{\mu_i}\leq \dim \oplus_{i=1}^{n} E_i.\label{proobraz subset E_1 E_2}
\end{equation}
Inequality (\ref{proobraz subset E_1 E_2}) implies
$$\sum_{i=1}^{n}\dim F_{\mu_i}\leq \sum_{i=1}^{n} \dim E_i.$$

Put in (\ref{S-D obsh beta}) $(u_1,u_2,\ldots,u_n)\in N$ . Then the left-hand side of equation (\ref{S-D obsh beta}) is equal to 1 and we have
\begin{equation}
1=\prod_{i=1}^{n}\prod_{j=1}^{n}\hat{\mu}_i(\beta_{ij}u_j),\quad (u_1,u_2,\ldots,u_n)\in N.\label{1=S-D obsh}
\end{equation}
Fix $j$. Then for each $u\in E_j$ there is $(u_1,u_2,\ldots,u_n)\in N$ such that $u_j=u$.
From this, (\ref{1=S-D obsh}), and $0\leq \hat{\mu}_i(y)\leq 1,y\in Y,$ it follows that $\hat{\mu}_i(\beta_{ij}u)=1,u\in E_j$. Hence, the following inclusions hold
$$\beta_{ij}(E_j)\subset F_{\mu_i},\quad i,j=1,2,\ldots,n.$$

Finally, we infer that the conditions of Lemma~\ref{Lem 2} are satisfied. Therefore $F_{\mu_i}=F$, where $F$ is a subgroup of $Y$, and $\beta_{ij}(F)=F,i,j=1,2,\ldots,n$.
$\blacksquare$

\begin{cor}\label{cor 1}
Let $Y$ be a finite Abelian group. Let $\hat{\mu}_i(y),i=1,2,\ldots,n,$ $n\geq 2$, be nonnegative characteristic functions on $Y$ satisfying equation $(\ref{S-D obsh beta})$, where $\beta_{1j}=\beta_{i1}=I,i,j=1,2,\ldots,n$.
Then either $F_{\mu_i}=\{0\},i=1,2,\ldots,n,$ or $F_{\mu_i}\neq\{0\},i=1,2,\ldots,n,$ and there is a nonzero subgroup $H$ of $Y$ such that $H\subset(\cap_{i=1}^{n} F_{\mu_i})$ and $\beta_{ij}H=H,i,j=1,2,\ldots,n$.
\end{cor}

\textbf{Proof. }
Assume that $F_{\mu_k}=\{0\}$ for some $k$. Fix a prime number $p$ and consider $Y_{(p)}$. Since $Y_{(p)}$ is a characteristic subgroup, we can consider the restriction of equality (\ref{S-D obsh beta}) to $Y_{(p)}$.
Then $Y_{(p)} \cap F_{\mu_k}=\{0\}$. From this and Lemma~\ref{Lem 3} it follows that $Y_{(p)} \cap F_{\mu_i}=\{0\},i=1,2,\ldots,n$. It means that each $F_{\mu_i}$ does not contain elements of order $p$. Since $p$ is arbitrary, we obtain $F_{\mu_i}=\{0\},i=1,2,\ldots,n$.

Suppose that $F_{\mu_k}\neq\{0\}$ for all $k$. Then, in particular, $F_{\mu_1}\neq \{0\}$. This implies that $Y_{(p)}\cap F_{\mu_1} \neq \{0\}$ for some $p$.
It follows from Lemma~\ref{Lem 3} that the subgroups $Y_{(p)}\cap F_{\mu_i},i=1,2,\ldots,n,$ are nonzero, they coincide, and they are invariant with respect to $\beta_{ij},i,j=1,2,\ldots,n$. Put $H=Y_{(p)}\cap F_{\mu_i}$. Then $H$ is desired subgroup.
$\blacksquare$

Next lemma is crucial for the proof of Theorem~\ref{Th 1}.

\begin{lem}\label{Lem 4}
Let $\xi_i, i=1,2,\ldots,n$, $n\geq 2$, be independent random variables with values in a group $X$ and distributions $\mu_i$ such that $\hat{\mu}_i(y)\geq 0$.
Consider the linear forms $L_j=\sum_{i=1}^{n}\alpha_{ij}\xi_i,$
where $\alpha_{ij}\in Aut(X),\alpha_{1j}=\alpha_{i1}=I,i,j=1,2,\ldots,n$. Suppose that the following condition is satisfied:

(\texttt{A}) For some $k$ any proper subgroup of $X$ does not contain the support of $\mu_{k}$.

Then the independence of $L_j$ implies that $\mu_i=m_X,i=1,2,\ldots,n$.
\end{lem}
\textbf{Proof. }
By Lemma~\ref{Lem 1} it follows that the equality
\begin{equation}
\prod_{i=1}^{n}\hat{\mu}_i\left(\sum_{j=1}^{n}\tilde{\alpha}_{ij}u_j\right)=\prod_{i=1}^{n}\prod_{j=1}^{n}\hat{\mu}_i(\tilde{\alpha}_{ij}u_j),\quad u_j\in Y,\label{S-D obsh}
\end{equation}
holds.

From (\texttt{\emph{A}}) it follows that
\begin{equation}
F_{\mu_k}=\{0\}. \label{D_i=0}
\end{equation}

Let $\pi\colon Y^n \rightarrow Y^n$ be a homomorphism defined by the formula
$$\pi(u_1,u_2,\ldots,u_n)=\left(\sum_{j=1}^{n}\tilde{\alpha}_{1j}u_j,\sum_{j=1}^{n}\tilde{\alpha}_{2j}u_j,\ldots,\sum_{j=1}^{n}\tilde{\alpha}_{nj}u_j\right),$$
where $u_j \in Y$.We will show that $\pi \in Aut(Y^n)$. Assume the converse, i.e. $\pi\not\in Aut(Y^n)$. Since $Y^n$ is a finite group, we obtain $Ker\pi\neq\{0\}$. Put in (\ref{S-D obsh}) $(u_1,u_2,\ldots,u_n)\in Ker{\pi},(u_1,u_2,\ldots,u_n)\neq0$:
\begin{equation}
1=\prod_{i=1}^{n}\prod_{j=1}^{n}\hat{\mu}_i(\tilde{\alpha}_{ij}u_j).\label{1=S-D obsh L4}
\end{equation}
From (\ref{1=S-D obsh L4}) and $\hat{\mu}_i(y)\geq 0$ it follows that all factors in the right-hand side of equation (\ref{1=S-D obsh L4}) are equal to 1. In particular, since $u_{j_0}\neq 0$ for some $j_0$, we obtain that $\hat{\mu}_i(\alpha_{ij_0}u_{j_0})=1,i=1,2,\ldots,n$, whence it follows that $F_{\mu_i}\neq \{0\},i=1,2,\ldots,n$. This contradicts condition (\ref{D_i=0}). Therefore, $\pi \in Aut(Y^n)$.

Let us prove that $\hat{\mu}_i(y)=0,i=1,2,\ldots,n,$ for all $y\in Y,y\neq 0$. Assume the converse. Then for some $l$ there is $\tilde{y}\neq 0$ such that
\begin{equation}
\hat{\mu}_l(\tilde{y})\neq 0 \label{mu_i neq 0}.
\end{equation}
Without loss of generality we can assume that $l=1$.

Putting in (\ref{S-D obsh}) $(\tilde{u}_1,\tilde{u}_2\ldots,\tilde{u}_n)=\pi^{-1}(\tilde{y},0,\ldots,0)$, we obtain:
\begin{equation}
\hat{\mu}_1(\tilde{y})=\prod_{i=1}^{n}\prod_{j=1}^{n}\hat{\mu}_i(\tilde{\alpha}_{ij}\tilde{u}_j).\label{predst 1}
\end{equation}
We note that there are at least two numbers $j_1,j_2$ such that $\tilde{u}_{j_1}\neq 0,\tilde{u}_{j_2}\neq 0$. Indeed, if $\tilde{u}_j=0,j=1,2,\ldots,n$, then we have the contradiction with $\pi^{-1}\in Aut(Y^n)$. If $\tilde{u}_{j_0}\neq 0,\tilde{u}_j=0,j\neq j_0,$ for some $j_0$, then $\pi(0,0,\ldots,\tilde{u}_{j_0},\ldots,0)=(\tilde{\alpha}_{1j_0}\tilde{u}_{j_0},\tilde{\alpha}_{2j_0}\tilde{u}_{j_0},\ldots,\tilde{\alpha}_{nj_0}\tilde{u}_{j_0})=(\tilde{y},0,\ldots,0)$. This contradicts $\tilde{\alpha}_{ij_0}\in Aut(Y)$. Hence, $\tilde{u}_{j_1},\tilde{u}_{j_2}\neq 0$ for some $j_1$ and $j_2$.
From inequalities
\begin{equation}
0\leq \hat{\mu}_i(y)\leq 1,\quad i=1,2,\ldots,n,\label{0 <mu_i< 1}
\end{equation}
and equation (\ref{predst 1}) we obtain
\begin{equation}
\hat{\mu}_1(\tilde{y})\leq\prod_{i=1}^{n}\hat{\mu}_i(\tilde{\alpha}_{ij_1}\tilde{u}_{j_1})\hat{\mu}_i(\tilde{\alpha}_{ij_2}\tilde{u}_{j_2}).\label{predst 1 ner}
\end{equation}
Put
\begin{equation}
C=\max_{1\leq i\leq n}\max_{y\neq 0}\hat{\mu}_i(y).\label{def C}
\end{equation}
By Corollary~\ref{cor 1} from (\ref{D_i=0}) we get
\begin{equation}
F_{\mu_i}=\{0\},\quad i=1,2,\ldots,n. \label{D_i=0 2}
\end{equation}
Combining (\ref{0 <mu_i< 1}), (\ref{mu_i neq 0}), and (\ref{D_i=0 2}), we obtain  $0<C<1$. Since $\tilde{u}_{j_1}\neq 0,\tilde{u}_{j_2}\neq 0$ and $\tilde{\alpha}_{ij_1},\tilde{\alpha}_{ij_2}\in Aut(Y),$ we have $\tilde{\alpha}_{ij_1}\tilde{u}_{j_1}\neq 0,\tilde{\alpha}_{ij_2}\tilde{u}_{j_2}\neq 0$. Hence, from (\ref{predst 1 ner}) and (\ref{def C}) we obtain that
$$\hat{\mu}_1(\tilde{y})\leq C^{2n}.$$

From (\ref{predst 1 ner}) and $\hat{\mu}_1(\tilde{y})\neq 0$ it follows that
\begin{equation}
\hat{\mu}_i(\tilde{\alpha}_{ij_1}\tilde{u}_{j_1}),\hat{\mu}_i(\tilde{\alpha}_{ij_2}\tilde{u}_{j_2})\neq 0,\label{dva somnozhitelya}
\end{equation}
where $\tilde{u}_{j_1}\neq 0,\tilde{u}_{j_2}\neq 0, i=1,2,\ldots,n$.

Using (\ref{dva somnozhitelya}) in the same way as (\ref{predst 1 ner}) was obtained from (\ref{mu_i neq 0}) we get an estimate for every factor in the right-hand side of (\ref{predst 1 ner})  and put this estimate in (\ref{predst 1 ner}). Repeating this process $m$ times we arrive at inequality that implies
$$\hat{\mu}_1(\tilde{y})\leq C^{{(2n)}^{m+1}}.$$
Since $C^{(2n)^{m+1}}\rightarrow 0$ as $m\rightarrow \infty$, we obtain $\hat{\mu}_1(\tilde{y})=0$. This contradicts the assumption.
Hence, $\hat{\mu}_i(y)=0,i=1,2,\ldots,n,$ for all $y\in Y,y\neq 0$. From this and (\ref{HF Haara}) we obtain that $\hat{\mu}_i(y)=\hat{m}_X(y), y\in Y,i=1,2,\ldots,n.$
Therefore, we have $\mu_i=m_X,i=1,2,\ldots,n$.
$\blacksquare$

\section{The proof of the main theorems}

\textbf{Proof of Theorem~\ref{Th 1}.}
Let $\delta_j\in Aut(X),j=1,2,\ldots,n$. Note that the linear forms $L_j=\sum_{i=1}^{n}\alpha_{ij}\xi_i,j=1,2,\ldots,n,$ are independent if and only if the linear forms $\delta_j L_j,j=1,2,\ldots,n,$ are independent. Since
$$L_j=\alpha_{1j}(\xi_1+\alpha_{1j}^{-1}\alpha_{2j}\xi_2+\ldots+\alpha_{1j}^{-1}\alpha_{nj}\xi_n),\quad j=1,2,\ldots,n,$$
without loss of generality we can assume that $\alpha_{1j}=I,j=1,2\ldots,n$, i.e.
\begin{equation}
L_j=\xi_1+\alpha_{2j}\xi_2+\ldots+\alpha_{nj}\xi_n,\quad j=1,2,\ldots,n.\label{formi}
\end{equation}
Put $\eta_i=\alpha_{i1}\xi_i$ and $\gamma_{ij}=\alpha_{ij}\alpha_{i1}^{-1}$. Then we can rewrite (\ref{formi}) in the form
$$L_1=\eta_1+\eta_2+\ldots+\eta_n,$$
$$L_j=\eta_1+\gamma_{2j}\eta_2+\ldots+\gamma_{nj}\eta_n,,\quad j=2,\ldots,n,$$
where random variables $\eta_i$ are independent.
Obviously, it suffices to prove Theorem~\ref{Th 1} assuming that
$\alpha_{1j}=\alpha_{i1}=I,i,j=1,2,\ldots,n$.

By Lemma~\ref{Lem 1} the functions $\hat{\mu}_i(y)$ satisfy equation (\ref{S-D obsh}).
Put $\nu_i=\mu_i\ast\bar{\mu_i},i=1,2,\ldots,n$. Then $\hat{\nu_i}(y)=|\hat{\mu}_i(y)|^2, y\in Y$. The functions $\hat{\nu_i}(y)$ are nonnegative and also satisfy equation (\ref{S-D obsh}).
We will prove that $\nu_i=m_K$, where $K$ is a subgroup of $X$. It is easy to see that this implies that $\mu_i=E_{x_i}\ast m_K,x_i\in X$, $i=1,2,\ldots,n$, i.e. $\mu_i\in I(X),i=1,2,\ldots,n$.

Put $F=\cap_{i=1}^n F_{\mu_i}$. Consider the set of subgroups $\{G_l\}\subset F$ such that $\tilde{\alpha}_{ij}G_l=\tilde{\alpha}_{ij},i,j=1,2,\ldots,n$. Denote by $H$ a subgroup of $Y$ such that $H$ is generated by all $\{G_l\}$.
It is not hard to prove that $H$ is a maximal subgroup of $Y$, which satisfies the condition

\texttt{(B)} $\hat{\nu_i}(y)=1,y\in \tilde{H},i=1,2,\ldots,n,$ $\tilde{\alpha}_{ij}\tilde{H}=\tilde{H},i,j=1,2,\ldots,n$.

Taking into account that $\hat{\nu}_i(y+h)=\hat{\nu}_i(y),i=1,2,\ldots,n,$ for all $y\in Y,h\in H,$ and the restrictions of the automorphisms $\tilde{\alpha}_{ij}$ of $Y$ to $H$ are automorphisms of $H$, consider the equation induced by equation (\ref{S-D obsh}) on the factor-group $Y/H$ putting $\tilde{\nu_i}([y])=\hat{\nu}_i(y),i=1,2,\ldots,n,$ and $\hat{\alpha}_{ij}[y]=[\tilde{\alpha}_{ij} y], y\in [y],[y]\in Y/H.$ Let $K=A(X,H)$. Note that $Y/H=(K)^{*}$. Thus, if we prove that $\tilde{\nu_i}([y])=\hat{m}_{K}([y]),[y]\in Y/H$, then we will obtain $\hat{\nu}_i(y)=\hat{m}_K(y), y\in Y,i=1,2,\ldots,n$.

Since $H$ is a maximal subgroup of $Y$, which satisfies condition \texttt{(B)}, we obtain that \{0\} is a maximal subgroup of $Y/H$, which satisfies condition \texttt{(B)} for the induced characteristic functions $\tilde{\nu}_i([y])$ and the induced automorphisms $\hat{\alpha}_{ij}$.

Hence, without loss of generality we suppose that
\begin{equation}
H=\{0\}. \label{U=0}
\end{equation}

Let us show that for some $k$ any proper subgroup of $X$ does not contain $\sigma(\nu_{k})$. This condition is equal to the condition $F_{\nu_k}=\{0\}$. Assume the converse. Then by Corollary~\ref{cor 1} there is a nonzero subgroup $\tilde{H}$ of the group $Y$ that satisfies condition \texttt{(B)}. This contradicts (\ref{U=0}). Hence, any proper subgroup of $X$ does not contain the support of $\nu_k$. Then by Lemma~\ref{Lem 4} $\nu_i=m_X,i=1,2,\ldots,n.$
$\blacksquare$

\begin{remark}\label{Remark 1}
The independence of the linear forms $L_j,j=1,2,\ldots,n,n\geq 2,$ where $\alpha_{1j}=\alpha_{i1}=I$ implies that $\xi_i=m_K\ast E_{x_i},i=1,2,\ldots,n$. Here, in contrast with the general case, the distributions of the random variables $\xi_i$ are the shifts of the Haar distribution on the same subgroup of $X$.
\end{remark}

Let us prove that Theorem~\ref{Th 1} is sharp in the following sense: in the class of finite Abelian groups the independence of $k$ linear forms of $n$ random variables, where $k<n$, does not imply that $\mu_i\in I(X)$.

\begin{thrm}\label{Th 2}
Let $n$ and $k$ satisfy the condition $n>k>1$. Let $X=(\mathbb{Z}(p))^n$, where $p>2$ is a prime number, such that $p$ does not divide $n$. Then there exist independent random variables $\xi_i, i=1,2,\ldots,n,$ with values in a group $X$ and distributions $\mu_i\not\in I(X),$ and automorphisms $\alpha_{ij}\in Aut(X),$ such that the linear forms $L_j=\sum_{i=1}^{n}\alpha_{ij}\xi_i,j=1,2,\ldots,k,$ are independent.
\end{thrm}

\textbf{Proof. }
It is obvious that it suffices to prove the statement for $k=n-1$.

Let $\alpha_{i,i-1}x=2x,x\in X,i=2,3,\ldots n$, and $\alpha_{ij}=I$ in other cases, $i=1,2,\ldots,n,j=1,2,\ldots,n-1$. It is clear that $\alpha_{ij}\in Aut(X)$.
Note that $Y\cong(\mathbb{Z}(p))^n,\tilde{\alpha}_{ij}=\alpha_{ij}$.

Let $e_1=(1,0,\ldots,0),e_2=(0,1,\ldots,0),\ldots,e_n=(0,0,\ldots,n)\in Y$.
Consider on $X$ the function
$$\rho_i(x)=1+Re(x,e_i).$$
Then  $\rho_i(x)\geq0,x\in X,$ and
$$\sum_{x\in X} \rho_i(x)m_X(\{x\})=1.$$
Denote by $\mu_i$ the distribution on $X$ with the density $\rho_i(x)$ with respect to $m_X$. We see that
\[
\hat{\mu}_i(y)=
\begin{cases}
1, & \text{$y$=0};\\
\frac{1}{2},& \text{$y=\pm e_i$;}\\
0, & y\in Y,y\not\in \{0,\pm e_i\}.
\end{cases}
\]
Obviously, $\mu_i \not\in I(X)$.
Let $\xi_i,i=1,2,..,n,$ be independent random variables  with values in $X$ and distributions $\mu_i$.
Let us show that the linear forms $L_j=\sum_{i=1}^{n}\alpha_{ij}\xi_i$ are independent. By Lemma~\ref{Lem 1} it suffices to show that the characteristic functions $\hat{\mu}_i(y)$ satisfy equation (\ref{S-D obsh nk}), which takes the form
$$\hat{\mu}_1(u_1+u_2+\cdots+u_{n-1})\hat{\mu}_2(2u_1+u_2+\cdots+u_{n-1})\cdots\hat{\mu}_n(u_1+u_2+\cdots+2u_{n-1})=$$
\begin{equation}
=\hat{\mu}_1(u_1)\hat{\mu}_1(u_2)\cdots\hat{\mu}_1(u_{n-1})\hat{\mu}_2(2u_1)\hat{\mu}_2(u_2)\cdots\hat{\mu}_2(u_{n-1})\cdots\label{S-D kontrpr}
\end{equation}
$$\cdots \hat{\mu}_n(u_1)\hat{\mu}_n(u_2)\cdots\hat{\mu}_n(2u_{n-1}).$$

Let us prove that the left-hand side of equation (\ref{S-D kontrpr}) does not equal to 0 if and only if $u_j=0,j=1,2,\ldots,n-1$.
Indeed, suppose that the left-hand side of (\ref{S-D kontrpr}) does not equal to 0. Then $u_j$ satisfy the system of equations
\begin{equation}\label{syst kontrpr}
\begin{cases}
u_1+u_2+\cdots+u_{n-1}=b_1,\\
2u_1+u_2+\cdots+u_{n-1}=b_2,\\
\ldots\ldots\ldots\\
 u_1+u_2+\cdots+2u_{n-1}=b_n,
\end{cases}
\end{equation}
where $b_i\in\{0,\pm e_i\}$.

From (\ref{syst kontrpr}) it follows that:

\begin{equation}\label{syst kontrpr resh}
\begin{cases}
\sum_{i=2}^{n}b_i=nb_1, \\
u_1=b_2-b_1,\\
u_2=b_3-b_1,\\
\ldots\ldots\ldots\\
u_{n-1}=b_n-b_1.
\end{cases}
\end{equation}

First equation of system (\ref{syst kontrpr resh}) implies that $b_i=0,i=1,2,\ldots,n$. Thus the unique solution of system (\ref{syst kontrpr}) is $u_j=0,j=1,2,\ldots,n-1$.

Taking into account that $\hat{\mu}_i(\pm e_j)=0$ for $i\neq j$, it easy to see that if  $u_j\neq 0$ for some $j$, then the right-hand side of equation (\ref{S-D kontrpr}) is equal to 0, i.e. the right-hand side of equation (\ref{S-D kontrpr}) does not equal to 0 if and only if $u_j=0,j=1,2,\ldots,n-1$. Hence, equality (\ref{S-D kontrpr}) holds for all $u_j\in Y$.
$\blacksquare$

Note that Theorem~\ref{Th 2} can be strengthened for $n=3$.
Denote by $G$ a group of the form
$(\ref{Lestnica})$. The following statements hold (\cite{Gr F}):

$1)$ Let $\alpha_i,\beta_i \in Aut(G),i=1,2,3$, $\xi_i$ be independent random variables with values in a group $X$ and distributions $\mu_i$. Suppose that linear forms $L_1=\alpha_1\xi_1+\alpha_2\xi_2+\alpha_3\xi_3$ and $L_2=\beta_1\xi_1+\beta_2\xi_2+\beta_3\xi_3$ are independent. If $X=G$, then all $\mu_i$ are degenerate distributions. If $X=\mathbb{Z}(3)\times G$, then either all $\mu_i$ are degenerate distributions or $\mu_{i_1}\ast E_{x_1}=\mu_{i_2}\ast E_{x_2}=m_{\mathbb{Z}(3)},x_i\in X,$ at least for two distributions $\mu_{i_1}$ and $\mu_{i_2}$.  If $X=\mathbb{Z}(5)\times G$, then either all $\mu_i$ are degenerate distributions or $\mu_{i_1}\ast E_{x_1}=m_{\mathbb{Z}(5)},x_1\in X,$ at least for one distribution $\mu_{i_1}$.

$2)$ If a group $X$ is not isomorphic to any of the groups mentioned in $1)$, then there exist $\alpha_i,\beta_i\in Aut(X),i=1,2,3$, and independent identically distributed random variables $\xi_i$ with values in $X$ and distribution $\mu \not\in I(X)$, such that the linear forms $L_1=\alpha_1\xi_1+\alpha_2\xi_2+\alpha_3\xi_3$ and $L_2=\beta_1\xi_1+\beta_2\xi_2+\beta_3\xi_3$ are independent.

\smallskip

We prove now that Theorem~\ref{Th 1} fails if $\alpha_{ij}$ are endomorphisms of $X$ and not all $\alpha_{ij}$ are automorphisms.

\begin{prop}\label{prop 1}
Assume that a group $X$ is not isomorphic to the group $\mathbb{Z}(p)$, where $p$ is a prime number. Then there are independent identically distributed random variables $\xi_1,\xi_2,$ with values in $X$ and distribution $\mu$ and nonzero endomorphisms $\alpha,\beta$ of $Y$ such that:

a) the linear forms $L_1=\alpha\xi_1+\beta\xi_2$ and $L_2=\xi_1+\alpha\xi_2$ are independent;

b) $\mu\not\in I(X)$;

c) $\sigma(\mu)=X$.
\end{prop}

\textbf{Proof.}
First we will show that there exist endomorphisms $\alpha,\beta$ of $X$ satisfying the conditions

1) $\alpha\not\in Aut(X),\beta \in Aut(X)$;

2) $\beta (Ker\alpha)=Ker\alpha$;

3) $\alpha^2x\neq \beta x$ for all $x\in X,x\neq 0$.

Without loss of generality we can suppose that $X$ is a $p$-primary group. By the structure theorem for finite Abelian groups
$$X=\prod_{k=1}^m(\mathbb{Z}(p^k))^{k_l},$$
where $k_l\geq0$. There are two possibilities: $X\cong\mathbb{Z}(p^k)$ and $X\not\cong\mathbb{Z}(p^k)$. If $X\cong\mathbb{Z}(p^k)$, where $k>1$, then put $\alpha x = px,x\in X$, $\beta=(p-1)x,x\in X$ . It easily can be proved that $\alpha$ and $\beta$ satisfy conditions 1)-3).

If $X\not\cong \mathbb{Z}(p^k)$, then $X=X_1\times X_2$, where $X_1,X_2$ are non-trivial subgroups of $X$. Denote by $(x_1,x_2),x_i\in X_i,$ elements of $X$. Put $\alpha (x_1,x_2)=(0,x_1),x\in X,\beta=I$. It is no hard to prove that conditions 1)-3) satisfied.

So let $\alpha$ and $\beta$ satisfy conditions 1)-3). It easily can be showed that a homomorphism $\pi: Y^2\rightarrow Y^2$ defined by the formula
\begin{equation}
\pi(u,v)=(\tilde{\alpha}u+v,\tilde{\beta} u+\tilde{\alpha}v) \label{def pi l5}
\end{equation}
is an automorphism of $Y^2$.
It is clear that $H=Ker\tilde{\alpha}\neq\{0\}$.
From (\ref{def pi l5}) and condition 2) it follows that $\pi H^2\subset H^2$. Since $\pi \in Aut(Y^2)$ and $Y^2$ is a finite group, we obtain
\begin{equation}
\pi H^2 =H^2.\label{piH=H l5}
\end{equation}

Put $K=A(X,H)$, $\mu=(1-b)m_X+bm_K$, where $0<b<1$. Then
\begin{equation}\label{def mu}
\hat{\mu}(y)=
\begin{cases}
1,& y=0,\\
b,& y\in H,y\neq 0,\\
0,& y\not\in H.
\end{cases}
\end{equation}

It is obvious that $\mu \not\in I(X)$ and $\sigma(\mu)=X$.

Consider independent identically distributed random variables $\xi_i,\xi_2$ with values in a group $X $ and with the distribution $\mu$.
We shall prove that $L_1$ and $L_2$ are independent. By Lemma~\ref{Lem 1} it suffices to show that the characteristic function $\hat{\mu}(y)$ satisfies equations (\ref{S-D obsh}) which takes the form
\begin{equation}
\hat{\mu}(\tilde{\alpha}u+v)\hat{\mu}(\tilde{\beta} u+\tilde{\alpha}v)=\hat{\mu}(\tilde{\alpha}u)\hat{\mu}(v)\hat{\mu}(\tilde{\beta} u)\hat{\mu}(\tilde{\alpha}v),\quad u,v\in Y. \label{S-D l5}
\end{equation}

If $u,v\in H$, then it is clear that (\ref{S-D l5}) holds.

We will show that if either $u\not\in H$ or $v\not\in H$ both sides of (\ref{S-D l5}) are equal to 0.

If either $u\not\in H$ or $v\not\in H$, then (\ref{def mu}) implies that the right-hand side of (\ref{S-D l5}) is equal to 0.
Let us show that the same is true for left-hand side of (\ref{S-D l5}).
Assume the converse. Then the following inclusions hold

\begin{equation}\label{vkl l5}
\begin{cases}
\tilde{\alpha}u+v\in H,\\
\tilde{\beta} u+\tilde{\alpha}v\in H.
\end{cases}
\end{equation}
The inclusions (\ref{vkl l5}) mean that $\pi(u,v)\in H^2$. Then (\ref{piH=H l5}) implies that $(u,v)\in H^2$, i.e. $u,v\in H$. This contradicts the assumption.
$\blacksquare$

The author would like to thank G.M.Feldman for the suggestion of the problems to me and useful comments and A.I.Illinsky for useful discussions and comments.

Mathematical Division\\
B. Verkin Institute for Low \\
Temperature Physics and Engineering\\
 of the National Academy \\
 of Sciences of Ukraine\\
47, Lenin Ave, Kharkov\\
61103, Ukraine

\end{document}